\documentclass{article}

\usepackage[preprint]{neurips_2023}

\usepackage{amsmath}
\usepackage{amssymb}
\usepackage{mathtools}
\usepackage{amsthm}
\usepackage[utf8]{inputenc} 
\usepackage[T1]{fontenc}    
\usepackage{hyperref}       
\usepackage{url}            
\usepackage{booktabs}       
\usepackage{amsfonts}       
\usepackage{nicefrac}       
\usepackage{microtype}      
\usepackage{xcolor}         
\usepackage{subfigure}

\theoremstyle{plain}
\newtheorem{theorem}{Theorem}[section]

\theoremstyle{definition}
\newtheorem{definition}[theorem]{Definition}

\usepackage[textsize=tiny]{todonotes}
\usepackage[textsize=tiny]{todonotes}

\title{High Dimensional Geometry and Limitations in System Identification}

\author{%
  Muhammad Abdullah Naeem and Miroslav Pajic\thanks{Muhammad Abdullah Naeem and Miroslav Pajic are with the Department of Electrical and Computer Engineering, 
        Duke University, Durham, NC 27708, USA, Email: 
        {\tt\small muhammad.abdullah.naeem@duke.edu, miroslav.pajic@duke.edu}}
}

\begin{document}

\maketitle
\begin{abstract}

We study the problem of identification of linear dynamical system from a single  trajectory, via excitations of isotropic Gaussian. In stark contrast with previously reported results, Ordinary Least Squares (OLS) estimator for even \emph{stable} dynamical system contains non-vanishing error in \emph{high dimensions}; which stems from the fact that realizations of non-diagonalizable dynamics can have strong \emph{spatial correlations} and a variance, of order $O(e^{n})$, where $n$ is the dimension of the underlying state space. Employing \emph{concentration of measure phenomenon}, in particular tensorization of \emph{Talagrands inequality} for random dynamical systems we show that observed trajectory of dynamical system of length-$N$ can have a variance of order $O(e^{nN})$. Consequently, showing some or most of the $n$ distances between an $N-$ dimensional random vector and an $(n-1)$ dimensional hyperplane in $\mathbb{R}^{N}$ can be close to zero with positive probability and these estimates become stronger in high dimensions and more iterations via \emph{Isoperimetry}. \emph{Negative second moment identity}, along with distance estimates give a control on all the singular values of \emph{Random matrix} of data, revealing limitations of OLS for stable non-diagonalizable and explosive diagonalizable systems.             
\end{abstract}
\section{Introduction}
The task of estimating parameters associated to an unknown dynamical systems is of fundamental importance in finance, econometrics, controls and reinforcement learning. Last decade has seen a tremendous surge in research work on finite-time analysis on learning for dynamical systems. In this work we demystify the working of OLS estimator for dynamical system $x_{k+1}=Ax_{k}+w_{k}$, i.e., estimating state transition matrix $A$ by using excitations $w_{k}$ of isotropic Gaussian and observed trajectory $(x_{k})_{k=1}^{N}$. Understanding this problem is causal for dealing with more involved problems like learning dynamics from input-output data (\cite{oymak2021revisiting}) and online Kalman filtering (\cite{tsiamis2022online}). In contrast with existing results \cite{sarkar2019finite,simchowitz2018learning, faradonbeh2018finite,faradonbeh2018finite} which only seems to work via careful selection of spectrum of state transition matrix and length of simulated trajectory and scalar systems, our analysis is particularly focused towards high dimensional state space. To the best of authors' knowledge we give first geometric analysis for system identification by leveraging upon recent developments in random matrix theory literature \cite{tao2010random,tao2012topics,vu2014modern,rudelson2009smallest} and direct sum decomposition of original state transition matrix onto its' invariant subspaces for non-diagonalizable systems. Indeed these invariant subspaces correspond to Jordan forms but standard results in literature only focus on algebraic properties of Jordan forms, ignoring the hidden geometric content. Along with Talagrand's concentration inequality we are able to conclude, as opposed to standard beliefs, main issue in system identification and for stable systems is strong spatial correlations: that stem from existence of invariant subspaces of large dimensions, \cite{naeem2023learning}. Large invariant subspaces are a consequence of discrepancy between algebraic and geometric multiplicity of eigenvalues associated to state transition matrix, leading to recursive generation of basis that naturally add spatial correlations. 
\paragraph{Contributions and main results:}
In this paper we show that analysis based on mere magnitude of eigenvalues associated to state transition is misleading as non-hermitian matrices naturally contain important geometric content that can not be taken for granted. Our analysis is based in Frobenius norm which is sensitive to all the modes rather than the most dominant one, as small error in operator norm can be misleading. We have been wasteful with some of the concentration results, as our aim is to make things interpretable and accessible to general machine learning community. If the observed trajectory is $[x_1,x_2, \ldots,x_N]=:\mathbb{X} \in \mathbb{R}^{n \times N}$, the trick is to understand the geometry of its rows $y_1,y_2, \ldots,y_{n}$. Consider the row $y_j$ and the hyperplane $v_j$ spanned by all rows modulo $y_j$, error increases as soon as any of these $n$ distances become small. In the case of diagonalizable dynamics this is the case when their is an explosive mode; as it can have large variance and with decent probability be close to a high dimensional random subspace. Non-diagonalizable stable dynamics are naturally correlated in space, and in higher dimensions with increase in iterations their variance and largest singular value increases exponentially.
\paragraph{Structure of the paper:}
We begin section \ref{sec:Not and Prelim} with notations, followed by important results from concentration of measure phenomenon for random dynamical systems. In Section \ref{sec:sys-id-diag}, we formally define the OLS problem and provide a very' detailed analysis on Diagonalizable systems, as it will simplify understanding for non-hermitaian dynamics. Section \ref{sec:non-diag} begins with geometric analysis of non-diagonalizable matrices and limitations of OLS are verified via simulation results. We conclude our analysis in Section \ref{sec:conc}, with conclusions and directions on future work. 
\section{Concentration of measure phenomenon for random dynamical systems}
\label{sec:Not and Prelim}
\paragraph{Notations:}We use ${I}_{n}\in\mathbb{R}^{n\times n}$ to denote the $n$ dimensional identity matrix. 
For random variables $x$ and $y$, $corr(x,y)$ denote the correlation/covariance. 
 $\mathcal{S}^{p-1}:=\{x \in \mathbb{R}^{p}:\|x\|_2=1\}$, is the unit sphere in $\mathbb{R}^{p}$. Whereas,  $A^{T}$, $\rho(A)$,  $\|A\|$, $\|A\|_{F}$, $det(A)$, $tr(A)$ and $\sigma(A)$   represent the  transpose, spectral radius,  matrix 2-norm (operator norm) , Frobenius norm, determinant, trace and set of eigenvalues(spectrum) of $A \in \mathbb{R}^{n \times n}$ respectively. $O$, $\Theta$ and $\Omega$ are standard order relations. 
For a positive definite matrix $A$, largest and smallest eigenvaues are denoted by $\lambda_{max}(A)$ and $\lambda_{min}(A)$, respectively. Associated with every rectangular matrix $\mathbb{X} \in \mathbb{R}^{n \times N}$ are its' singular values $\sigma_{1}(\mathbb{X}) \geq \sigma_{2}(\mathbb{X}), \ldots, \sigma_{n}(\mathbb{X}) > 0$, where without loss of generality we assume that $N>n$. More frequently used are the largest singular value,
$\sigma_{1}(\mathbb{X}):= \sup_{a \in \mathcal{S}^{N-1}} \|\mathbb{X}a\| $ and the smallest singular value $\sigma_{n}(\mathbb{X}):= \inf_{a \in \mathcal{S}^{N-1}} \|\mathbb{X}a\|$. Condition number of a matrix $\mathbb{X}$ is the ratio of the largest and least singular value, denoted by $\kappa(\mathbb{X})= \frac{\sigma_{1}(\mathbb{X})}{\sigma_{n}(\mathbb{X})}$. If the span of image space of $\mathbb{X}$ is $\mathbb{R}^{n}$, more compactly written as $Im(\mathbb{X})=\mathbb{R}^{n}$, least singular value equals the inverse of the norm of inverse of matrix $\mathbb{X}$ i.e., $\sigma_{n}(\mathbb{X})= \frac{1}{\|\mathbb{X}^{-1}\|}$. A function $g: \mathbb{R}^{n} \rightarrow \mathbb{R}^{p}$ is Lipschitz with constant $L$ if for every $x,y \in \mathbb{R}^{n} $, $\|g(x)-g(y)\| \leq L\|x-y\|$. Space of probability measure on  $\mathcal{X}$(continuous space) is denoted by  $\mathcal{P(\mathcal{X})}$ and space of its Borel subsets is represented by $\mathbb{B}\big(\mathcal{P(\mathcal{X})}\big)$. For a function $f$ and $\mu \in \mathcal{P(X)}$, we use $<f>_{\mu}$ to denote expectation of $f$ w.r.t $\mu$.On a metric space $(\mathcal{X},d)$, for $\mu, \nu \in \mathcal{P(\mathcal{X})}$, we define Wasserstein metric of order $p \in [1, \infty)$~as
\begin{equation}
\label{eq:WM}
    \mathcal{W}_{d}^{p} (\nu,\mu)= \bigg(\inf_{(X,Y) \in \Gamma(\nu,\mu)} \mathbb{E}~d^{p}(X,Y)\bigg)^{\frac{1}{p}};
\end{equation}
here, $\Gamma(\nu,\mu) \in P(\mathcal{X}^{2})$, and $(X,Y) \in \Gamma(\nu,\mu)$ implies that random variables $(X,Y)$ follow some probability distributions on $P(\mathcal{X}^{2})$ with marginals $\nu$ and $\mu$. Another way of comparing two probability distributions on $\mathcal{X}$ is via relative entropy, which is defined as
\begin{equation}
\label{eq:ent} 
    Ent(v||u)=\left\{ \begin{array}{lr}
    \int \log\bigg(\frac{d\nu}{d\mu}\bigg) d\nu, & \text{if}~ \nu << \mu,
         \\ +\infty, & \text{otherwise}. 
         \end{array}\right.
\end{equation}
 
\paragraph{Notions of stability:}
Markov chain under consideration is $n$ dimensional LG with isotropic noise:
\begin{equation}
\label{eq:LGS_1}
    x_{t+1}= Ax_t+ w_{t}, \hspace{10pt} \rho(A) <1 \hspace{10pt} \text{and i.i.d }~ w_{t} \thicksim \mathcal{N}(0,\mathcal{I}_n).
\end{equation}
It mixes to stationary distribution $\mu_{\infty} \thicksim \mathcal{N}(0, P_{\infty})$, where the controllability grammian $P_{\infty}$ is the unique positive definite solution of the following Lyapunov equation:
\begin{equation}
\label{eq:contgram}
    A^{T}P_{\infty}A-P_{\infty}+I_{n}=0.
\end{equation}
Stability in controls community for Linear systems correspond to  $\rho(A)<1$ (marginally stable corresponds $\rho(A) \geq 1$ and explosive system when $\rho(A)>1$) 
\begin{definition} [Talagrands inequality]Consider metric space $(\mathcal{X},d)$ and reference  probability measure $\mu \in P(\mathcal{X})$. Then we say that $\mu$ satisfies  $\mathcal{T}_{1}^{d^2} (C)$ or to be concise $\mu \in  \mathcal{T}_{1}^{d^2} (C)$ for some $C>0$ if for 
all $\nu \in P(\mathcal{X})$ it holds~that 
\begin{equation}
    \label{eq:t1}
    \mathcal{W}_{d}^{2} (\mu, \nu) \leq \sqrt{2 C Ent(\nu||\mu)}.
\end{equation}
\end{definition}

\begin{theorem}
\label{lm:b-g}
If $\mu$ satisfies $\mathcal{T}_{1}^{d^2} (C)$ then for all Lipschitz function $f$ with $<f>_{\mu}:= \mathbb{E}_{\mu} f$, it holds that
\begin{flalign}
    \label{eq:bg} & \int e^{\lambda(f- <f>_{\mu})} d\mu \leq \exp(\frac{\lambda^2}{2}C \|f\|_{L(d)} ^2), \hspace{15pt} \text{where~~~~}  \|f\|_{L(d)}:= \sup_{x \neq y} \frac{|f(x)-f(y)|}{d(x,y)}.
\end{flalign}
Moreover, (\ref{eq:bg}) combined with the Markov inequality implies that if we sample $x$ from $\mu \in \mathcal{T}_{1}^{d^2}(C)$, then
\begin{equation}
    \label{eq:iid} \mathbb{P} \big[ | f(x) -<f>_{\mu} | > \epsilon \big] \leq 2 e^{\big(-\frac{ \epsilon ^2}{2C \|f\|_{L(d)}^2}\big)}.
\end{equation}
\end{theorem}
Since we will be dealing with the dependent processes, it is important to understand the concentration phenomenon for the entire length of simulated trajectory, a tensorization procedure for Talagrands inequality is required. If one is concerned with the concentration of time averages around steady state spatial average, as in \cite{naeem2022transportation} an $\ell_{1}$ based metric suffices. However, when dealing with projections and distances in large dimension, following metric on product space is more appropriate with associated variance constants.
\begin{theorem}
\label{thm:ten_tal}[Tensorization of Talagrands' Inequality for Dependent Processes]
Let, $d_{(N)} ^2 (x^N,y^N):= \sqrt{\sum\nolimits_{i=1}^{N} d^2(x_i,y_i)}$ be an $\ell^{2}$ inspired metric on $\mathcal{X}^{N}$.  
If system transition matrix, $A$ for LG satisfies $\|A\|<1$, then the process level law of $(x_1,x_2, \ldots,x_{N})$ satisfies $T_{1}^{d_{(N)}^2} \in O\bigg(\frac{1}{(1-\|A\|)^2}\bigg)$, if  $\|A\|=1$ we have $T_{1}^{d_{(N)}^2} \in O(N(N+1))$  and for $\|A\|>1$, $T_{1}^{d_{(N)}^2} \in O(\|A\|^{N} N)$ \cite{blower2005concentration}. 
\end{theorem}
A fundamental limitation of existing results, stem from the fact that they base their analysis only on the extreme singular values of the data matrix which do not offer much of geometric insight. This brings us to an extremely important, although simple to prove equality that relates the distance between different rows (whose concentration we understand well because of preceding theorem \ref{thm:ten_tal}) of the data matrix to all of its' singular values. 
\begin{theorem}
\cite{tao2010random}
[Negative second moment identity \label{thm:neg_2nd_moment_ineq}] Let $1 \leq d \leq p$ and $\mathbb{Y} \in \mathbb{R}^{d \times p}$ be a full rank matrix with singular values, $\sigma_{1}(\mathbb{Y}) \geq \sigma_{2}(\mathbb{Y}) \ldots \geq \sigma_{d}(\mathbb{Y}) $. Let $V_{j}$ be the hyperplane generated by all the rows of $\mathbb{Y}$-except the $j-th$ : i.e., span of $Y_1,Y_2, \ldots, Y_{j-1}, Y_{j+1}, \ldots, Y_{d}$  for $1 \leq j \leq d$ and $(e_{j})_{j=1}^{d}$ be the canonical basis of $\mathbb{R}^{d}$
 \begin{equation}
 \label{eq:neg_2nd_moment_ineq}
     \sum_{j=1}^{d} \frac{1}{\sigma_{j}^2(\mathbb{Y})}=\sum_{j=1}^{d}\big \langle \big(\mathbb{Y}\mathbb{Y}^{T}\big)^{-1}e_{j} ,e_{j} \big \rangle =\sum_{j=1}^{d} \frac{1}{d^{2}(Y_j, V_j)}. 
 \end{equation}
\end{theorem}
As the quantifying distance between a subspace and vector involves taking projections, high probability estimates will require concentration behavior of projections of isotropic Gaussians, as in: 
\begin{theorem}
\label{thm:gausssubcon}[Gaussian Projection and Isoperimetry, Lemma 3.2 in \cite{barvinok2005math}]
    Let $\gamma_{n}$ be isotropic Gaussian in $\mathbb{R}^{n}$  and  $ S \subset \mathbb{R}^{n}$ be a $k-$ dimensional subspace. Given $x \in \mathbb{R}^{n}$, let $x_{S}$ denote the projection of $x$ onto $S$. Then for any $\delta \in (0,1)$
\begin{align}
& \nonumber \gamma_{n}\bigg(x \in \mathbb{R}^{n} : \frac{\|x_{S}\|}{\|x\|} \geq (1-\delta)^{-1}\sqrt{\frac{k}{n}}\bigg) \leq e^{-\frac{\delta^2 k}{4}}+ e^{-\frac{\delta^2 n}{4}} \\ & \label{eq:concsubspace}\gamma_{n}\bigg(x \in \mathbb{R}^{n} : \frac{\|x_{S}\|}{\|x\|} \leq (1-\delta)\sqrt{\frac{k}{n}}\bigg) \leq e^{-\frac{\delta^2 k}{4}}+ e^{-\frac{\delta^2 n}{4}},
\end{align}
and we will often denote this behavior as $\frac{\|x_{S}\|}{\|x\|} \thicksim \sqrt{\frac{k}{n}}$ or say with overwhelming probability $\frac{\|x_{S}\|}{\|x\|}=\Theta(\sqrt{\frac{k}{n}})$  
\end{theorem}

\section{Problem statement: system Identification via single trajectory}
\label{sec:sys-id-diag}
\subsection{Ordinary Least Squares estimator}
In this section we analyse the problem of OLS estimation for system transition matrix $A$ from single observed (as in \cite{sarkar2019near}, \cite{simchowitz2018learning}, \cite{tsiamis2021linear}) trajectory of $(x_0,x_1, \ldots,x_{N})$ satisyfing:
\begin{equation}
    \label{eq:LGS}
    x_{t+1}=Ax_{t}+\eta_{t}, \hspace{10pt} \text{ where } \eta_{t} \thicksim \mathcal{N}(0,\mathcal{I}_n) .
\end{equation}  
Here $\eta_{t}$ acts as \emph{excitation} signal, so that dynamical system is well-explored. 
OLS solution is:
\begin{equation}
    \label{eq:OLSsol} \hat{A}= \arg \min_{B \in \mathbb{R}^{n \times n}}  \sum_{t=0}^{N-1} \|x_{t+1}-Bx_{t}\|.
\end{equation}
Let $\mathbb{X}_{+}=[x_1, x_2,  \ldots, x_N]$ and $ \mathbb{X}_{-}=[x_0, x_1, \ldots, x_{(N-1)}]$, noise covariates $E=[\eta_0, \eta_1, \ldots, \eta_{N-1}]$, $y_{j}$ be the rows of $\mathbb{X}_{-}$ and $v_{j}$ be the hyperplane as defined in theorem \ref{thm:neg_2nd_moment_ineq}. Then the closed form expression for  Least squares solution and error are:
\begin{align}
 & \label{eq:OLS}    \hat{A}= \mathbb{X}_{+}\mathbb{X}_{-} ^{T}(\mathbb{X}_{-}\mathbb{X}_{-}^{T})^{-1}  \\ & \label{eq:OLSerror} \|A-\hat{A}\|_{F} =\|E\mathbb{X}_{-}^{T} (\mathbb{X}_{-}\mathbb{X}_{-}^{T})^{-1}\|_{F}
\end{align}
\paragraph{Error bounds and statistics of Gaussian Ensemble:}
Spectral statistics of a random matrix of i.i.d Gaussian ensembles  $E$ (each entry of the matrix is normally distributed with mean $0$ and variance $1$), has been very well studied, and we would like to leverage upon that information and by definition $\|A-\hat{A}\|_{F}^{2}= \sum_{k=1}^{n} \sigma_{k}^{2}\big[E\mathbb{X}_{-}^{T} (\mathbb{X}_{-}\mathbb{X}_{-}^{T})^{-1}\big]$, using Courant-Fischer:
\begin{align}
    \sigma_{n}(E_{N(\mathbb{X}_{-})^{\perp}}) \sigma_{k}\big[\mathbb{X}_{-}^{T} (\mathbb{X}_{-}\mathbb{X}_{-}^{T})^{-1}\big] & \nonumber  \leq \sigma_{k}\big[E\mathbb{X}_{-}^{T} (\mathbb{X}_{-}\mathbb{X}_{-}^{T})^{-1}\big] := \max_{V \subset \mathbb{R}^{n}: dim(V)=k} \bigg(\min_{ x \in V \cap \mathcal{S}^{n-1} }\|E\mathbb{X}_{-}^{T} (\mathbb{X}_{-}\mathbb{X}_{-}^{T})^{-1}x\|\bigg)\\ & \nonumber=
    \min_{V \subset \mathbb{R}^{n}: dim(V)=n-k+1} \bigg(\max_{ x \in V \cap \mathcal{S}^{n-1}}\|E\mathbb{X}_{-}^{T} (\mathbb{X}_{-}\mathbb{X}_{-}^{T})^{-1}x\| \bigg) \\ & \nonumber \leq \sigma_{1}(E_{N(\mathbb{X}_{-})^{\perp}}) \sigma_{k}\big[\mathbb{X}_{-}^{T} (\mathbb{X}_{-}\mathbb{X}_{-}^{T})^{-1}\big]
\end{align}
where, owing to a very important observation by \cite{simchowitz2018learning}: $\mathbb{X}_{-}^{T} (\mathbb{X}_{-}\mathbb{X}_{-}^{T})^{-1}$ is a bijection from  $Im(\mathbb{X}_{-})= \mathbb{R}^{n}$ to $N(\mathbb{X}_{-})^{\perp}$ (subspace orthogonal to null space of $\mathbb{X}_{-}$) and $dim[N(\mathbb{X}_{-})^{\perp}]=n$, so $E_{N(\mathbb{X}_{-})^{\perp}}$ should be viewed as an $n \times n$ square matrix of i.i.d Gaussian ensembles. So now we are left with analyzing:
\begin{equation}
    \|\mathbb{X}_{-}^{T} (\mathbb{X}_{-}\mathbb{X}_{-}^{T})^{-1}\|_{F}^{2}=Tr[\mathbb{X}_{-}^{T} (\mathbb{X}_{-}\mathbb{X}_{-}^{T})^{-2}\mathbb{X}_{-}]=Tr[(\mathbb{X}_{-}\mathbb{X}_{-}^{T})^{-1}]= \sum_{j=1}^{n} \frac{1}{\sigma_{j}^2(\mathbb{X}_{-})},
\end{equation}
using negative second moment identity, we conclude:

\begin{align}
    \label{eq:ubderrlst} \sigma_{n}(E_{N(\mathbb{X}_{-})^{\perp}}) \sqrt{\sum_{j=1}^{n} \frac{1}{d^{2}(y_{j},v_j)} } \leq \|A-\hat{A}\|_{F}  \leq \sigma_{1}(E_{N(\mathbb{X}_{-})^{\perp}}) \sqrt{\sum_{j=1}^{n} \frac{1}{d^{2}(y_{j},v_j)}}.    
\end{align}

  It is well know that(see e.g., \cite{cook2017selected}), $\sigma_{1}(E_{N(\mathbb{X}_{-})^{\perp}}) \thicksim \sqrt{n}$  and  $\sigma_{n}(E_{N(\mathbb{X}_{-})^{\perp}}) \thicksim \frac{1}{\sqrt{n}}$. Therefore, we only need to study the concentration behavior of $\sqrt{\sum_{j=1}^{n} \frac{1}{d^{2}(y_{j},v_j)} }$, which brings us to the following section:  
\subsection{Learning rate for stable Diagonalizable dynamics}
\paragraph{No spatial correlations for Hermitian state transition matrix: } 
Recall that if $A=A^{T}$, then there exists an orthogonal matrix $U \in \mathbb{R}^{n \times n}$ and a diagonal matrix $\Omega \in \mathbb{R}^{n \times n}$ such that $A=U^{T} \Omega U$ and consequently a similarity transform for dynamical system given in equation \eqref{eq:LGS}, i.e., $x_{t} \mapsto Ux_t=z_t= [z_{t,1}, \ldots, z_{t,n}], \hspace{3pt} z_{t+1}=\Omega z_t+U \eta_t $, where $U\eta_t$ is again an isotropic Gaussian and we essentially have $n$, one dimensional LGs i.e., for every $t \in \mathbb{N}$ and $l \neq m \in [1,2, \ldots, n] $, $z_{t,l}$ and $z_{t,m}$ are not correlated to each other implying rows of $\mathbb{X}_{-}$ are also independent. Moreover if $\rho(A)<1$ (i.e., magnitude of eigenvalues of $A$: $\lambda_{m}$ is less than 1) then $\|A^{k}\|_{2}<1$ for all $k \in \mathbb{N}$, implying for each $s \in \mathbb{Z}_{+}$ and $m \in [1,2, \ldots, n]$, temporal correlation of $z_{t+s,m}$ and $z_{t,m}$, decays exponentially fast i.e., $corr(z_{t,m},z_{t+s,m})= O\big(\frac{1}{|\alpha_{m}|^s}\big)$, where $ \frac{1}{|\alpha_{m}|}=|\lambda_{m}|$ and  every row satisfies dimension free Talagrands' inequality(although dependent on spectral radius, $\rho$). Leveraging upon Talagrands-inequality we will show, distance-squared between each column $Y_{j}$ and subspace $V_j$ concentrates sharply around $\bigg[ (N-n+1)-\frac{\sqrt{N-n+1}}{1-\rho}, (N-n+1)+\frac{\sqrt{N-n+1}}{1-\rho} \bigg]$,  \emph{this should be interpreted as: probability of $d^2$ is atleast $\epsilon$ standard deviations $\big(\frac{\sqrt{N-n+1}}{1-\rho}\big)$ away from its' mean $(N-n+1)$ decays as $C\exp(-c \epsilon^2)$ for some $c,C>0$}  and consequently:

\begin{theorem}
 \label{thm:diagerror}   
 Suppose that $\rho:=\rho(A)<1$ and the dynamical system is diagonalizable, then with high probability
 
\begin{equation}
    \Omega\bigg(\sqrt{\frac{1-\rho}{(1-\rho)(N-n+1)+ \sqrt{N-n+1}}\bigg)} \leq \|A-\hat{A}\|_{F} \leq  O\bigg(\sqrt{\frac{(1-\rho)n^2}{(1-\rho)(N-n+1)-\sqrt{N-n+1}}}\bigg)  \end{equation}
\end{theorem}

\begin{proof}
    Our argument is along the lines of Corollary 2.1.19 of \cite{tao2012topics}, but with dependent covariates. Recall that $y_{j}$ is an $N$ dimensional Gaussian with covariance matrix $\Sigma_{N}$, so we can also write $y_{j}=\Sigma_{N}^{\frac{1}{2}}z_{N}$. where $z_{N}$ is $N-$ dimensional isotropic Gaussiann. Tensorization argument in theorem \ref{thm:ten_tal} restricted to one-dimensional stable dynamical system reveals $\|\Sigma_{N}^{\frac{1}{2}}\|_{2}=O\big(\frac{1}{1-\rho}\big)$. Simple, linear algebra reveals that $d^{2}(y_{j},v_{j})=  \langle P_{v_{j}^{\perp}}y_{j}, P_{v_{j}^{\perp}}y_{j} \rangle$, where $P_{v_{j}^{\perp}}$ is  orthogonal projection onto $v_{j}^{\perp}$ and $\mathbb{E}d^{2}(y_{j},v_{j})=Tr(\Sigma_{N}^{\frac{1}{2}}P_{v_{j}^{\perp}}\Sigma_{N}^{\frac{1}{2}})=\Theta(N-n+1)$. We know by Talagrands inequality:
\begin{equation}
    \mathbb{P}(|d(y_j,v_j) - \mathbb{E}d(y_j,v_j)|\geq \epsilon) \leq 2e^{-\frac{\epsilon^{2}(1-\rho)^{2}}{2}}.
\end{equation}
Moreover, by Gaussian projection lemma $\|P_{v_{j}^{\perp}}\Sigma_{N}^{\frac{1}{2}}z_{N}\|=O(\sqrt{N-n+1})$, with overwhelming probability:
\begin{equation}
    \mathbb{P}(|d^2(y_j,v_j) - \mathbb{E}d^2(y_j,v_j)|\geq \epsilon \sqrt{N-n+1}) \leq 2e^{-\frac{\epsilon^{2}(1-\rho)^{2}}{2}}.
\end{equation}
 So $d^2(y_j,v_j)$ concentrates in $  \big[ (N-n+1)-\frac{\sqrt{N-n+1}}{1-\rho}, (N-n+1)+\frac{\sqrt{N-n+1}}{1-\rho} \big]$, combined with $\sigma_{1}(E_{N(\mathbb{X}_{-})^{\perp}}) \thicksim \sqrt{n}$, $\sigma_{n}(E_{N(\mathbb{X}_{-})^{\perp}}) \thicksim \frac{1}{\sqrt{n}}$ and independence of the rows; result follows.
\end{proof}

\begin{figure}[!t]
\centering
\parbox{7.5cm}{
\includegraphics[width=7cm]{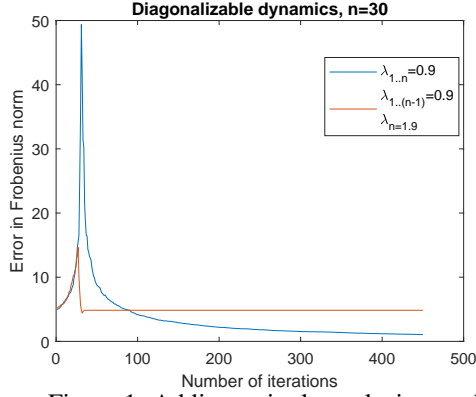}
\vspace{-10pt}
\caption{Adding a single explosive mode}
\label{fig:2figsA}}
\end{figure}

\paragraph{Statistics of singular values of data matrix:} \begin{theorem}
Given a Hermitian, stable state-transition matrix with spectral radius $\rho$, Largest singular value of the data matrix, $\sigma_{1}(\mathbb{X}_{-}) \thicksim \frac{\sqrt{N}+\sqrt{n}}{\sqrt{1-\rho^2}} $. There exists constants $c,C>0$ independent of the dimension of the state space or number of iterations such that for every $\delta>0$
\begin{align}
    \mathbb{P}\bigg(\sigma_{1}(\mathbb{X}_{-}) \geq \delta \frac{\sqrt{6}}{1-\rho} ( \sqrt{N}+\sqrt{n})\bigg) \leq \mathbb{P}\bigg(\sigma_{1}(\mathbb{X}_{-}) \geq \delta \sqrt{\frac{6}{1-\rho^2}} ( \sqrt{N}+\sqrt{n})\bigg) \leq C e^{-c\delta^2 ( \sqrt{N}+\sqrt{n})^2}
\end{align}
\begin{proof}
    Proof is an extension of standard concentration result for largest singular value of random rectangular matrix with i.i.d entries, see e.g., \cite{rudelson2009smallest}. One pays an entropy cost of $\bigg(\frac{3}{\epsilon}\bigg)^{N}$ to convert uncountable supremum to a finite one. Only notable extension is:
    $\mathbb{X}_{-}a$ for every $a \in S^{N-1}$ is an $n-$ dimensional Gaussian with spatially independent enteries and elementwise variance is upper bounded by $1+\rho^2+\rho^{4}+ \ldots \rho^{2(N-1)}$ that we will upper bound by $\frac{1}{1-\rho^2}$ so this is true for all $N$. For any $K>0$:

\begin{align}
\mathbb{P}\bigg(\sigma_{1}(\mathbb{X}_{-}) \geq K( \sqrt{N}+\sqrt{n})\bigg) & \label{eq:sup-enet} \leq \bigg( \frac{3}{\epsilon}\bigg)^{N} \inf_{s>0}\mathbb{P}\bigg(e^{s\|\mathbb{X}_{-}a\|_{2}} \geq e^{sK(1-\epsilon)^2( \sqrt{N}+\sqrt{n})^2}\bigg) \\ & \label{eq:utailnorm} \leq \frac{6^{N}3^{\frac{n}{2}}}{\exp\bigg( \frac{1-\rho^2}{3} K^{2} (\sqrt{N}+\sqrt{n})^2  \bigg)}, 
\end{align}    
where $a$ is a typical point in $\epsilon$ maximal net and \ref{eq:utailnorm} follows by upper tail concentration of $n$-dimensional Gaussian with independent elements each of variance $\frac{1}{1-\rho^2}$ and Markov Inequality where among feasible $s$ we picked $s=\frac{1-\rho^2}{3}$ and $\epsilon=\frac{1}{2}$ and proof follows. 
\end{proof}
\paragraph{Smallest singular value and bulk of the spectrum:}
Quantifying behavior of the smallest singular value of a general random matrix is a very involved problem, as it is sensitive to scaling of $N,n$ and tend to be extremely difficult for discrete distributions. However, as any sane algorithm would require number of iterations to be greater than state space dimension, results on this scaling are well know. Secondly, concentration of all the distances, implies singular values are evenly space and the least singular value is $\sigma_{n}(\mathbb{X}_{-}) \thicksim \sqrt{N}-\sqrt{n-1}$.
\end{theorem}
\paragraph{Adding explosive mode and deterioration of OLS in High Dimensions:} Tensorization of Talagrands' inequality, given in Theorem (\ref{thm:ten_tal}) implies that the distance-squared between norm-unstable trajectory $y_j  \in \mathbb{R}^{N}$ with covariance $\Sigma_{N}$ and $n-1$ dimensional subspace $v_j$ concentrates around $\big[Tr(\Sigma P_{v_{j}^{\perp}})-e^{\frac{N}{2}} \sqrt{N}\sqrt{Tr(\Sigma P_{v_{j}^{\perp}})}, Tr(\Sigma P_{v_{j}^{\perp}})+e^{\frac{N}{2}} \sqrt{N} \sqrt{Tr(\Sigma P_{v_{j}^{\perp}})} \big]$. Let $\lambda_1 \geq \lambda_2 \geq \ldots \geq \lambda_{N} \geq 0$ be the eigenvalues of $\Sigma$ in descending order, where $\lambda_1=O(e^{N}N)$. So if the sum of $(N-n+1)$ eigenvalues of the covariance matrix with corresponding eigenvectors in the span of $v_{j}^{\perp}$, is small then one of the distances is small with positive probability and OLS will fail to provide a good estimate. This is also shown in figure \ref{fig:2figsA}, with dimension of the state space fixed $n=30$ and first simulation with only eigenvalue of $\lambda=0.9$. In the second run we add one explosive mode of $1.9$ by removing one of the previously stable mode of $0.9$ and OLS performance deteriorates.

\section{Structure of Invariant Subspaces associated with Non-Hermitian Operator }
\label{sec:non-diag}
\begin{figure}
    \centering
    \subfigure[]{\includegraphics[width=0.4\textwidth]{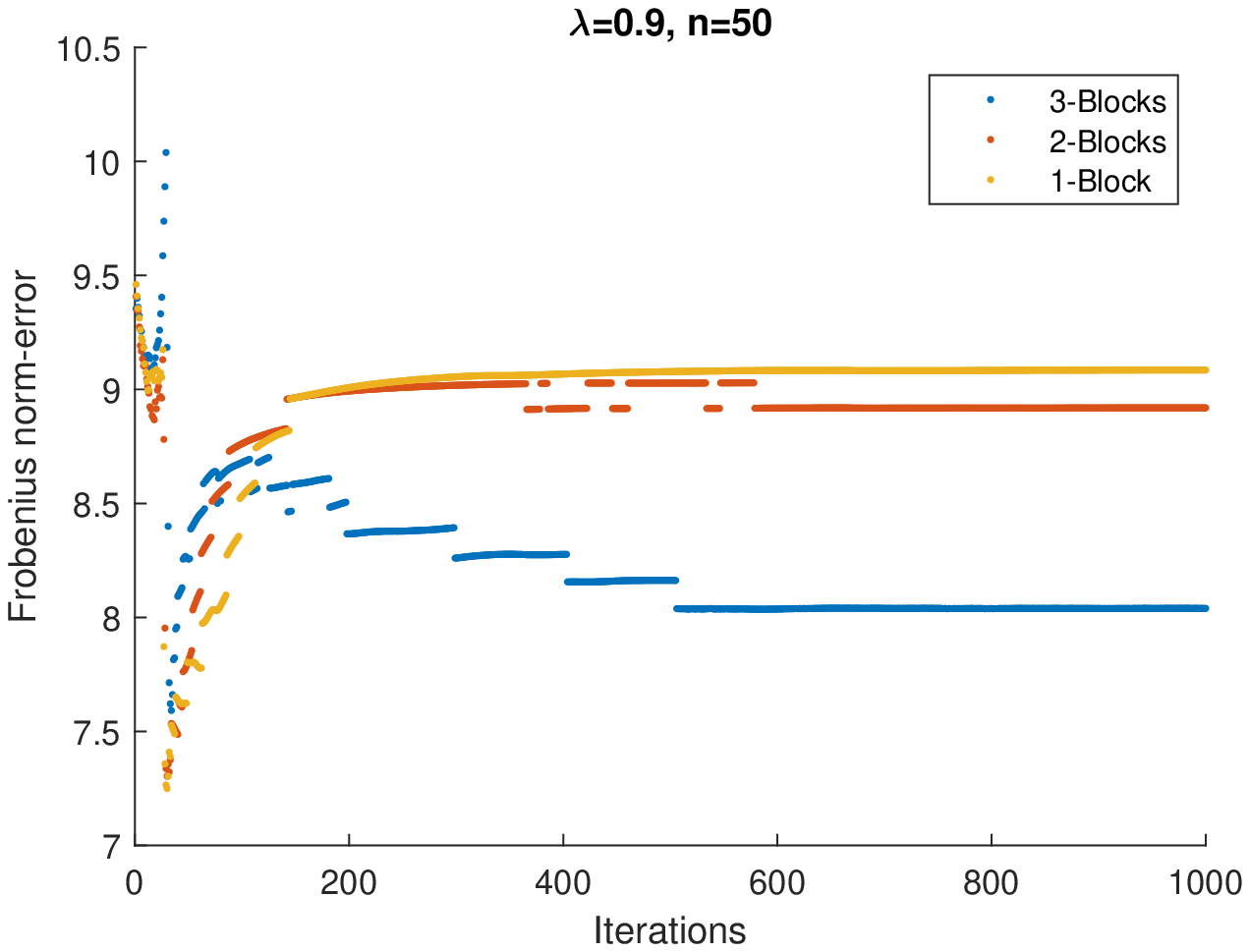}} 
    \subfigure[]{\includegraphics[width=0.4\textwidth]{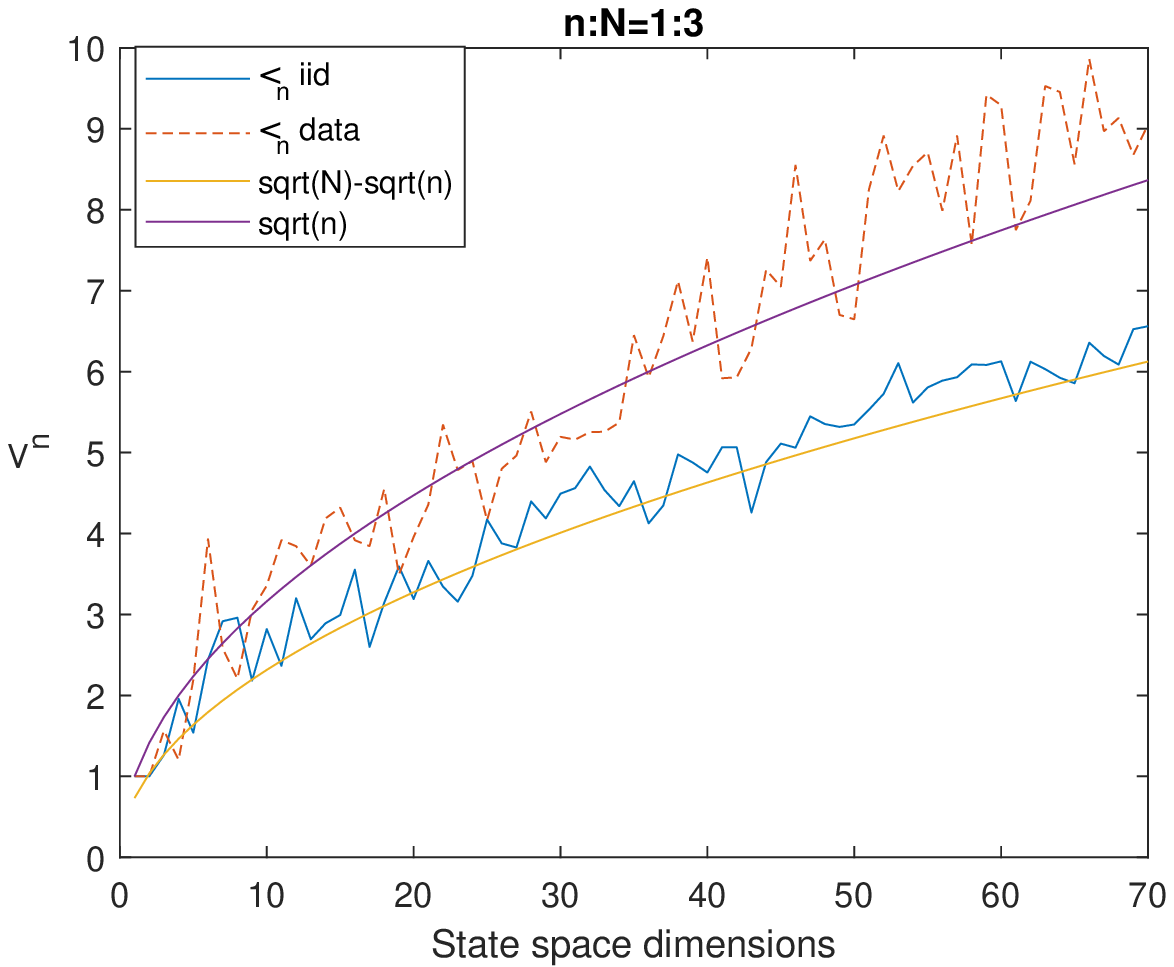}}
    \subfigure[]
    {\includegraphics[width=0.4\textwidth]
    {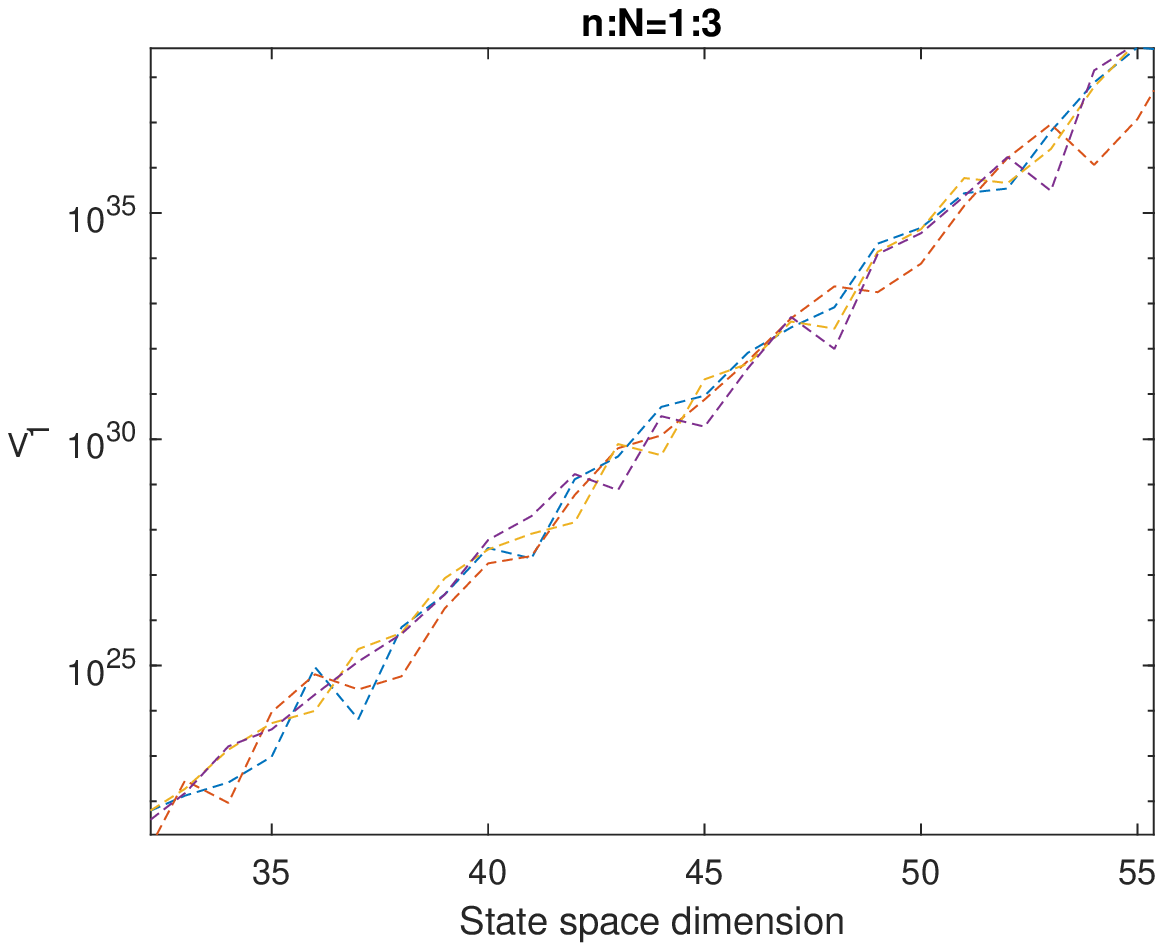}}
    \subfigure[]{\includegraphics[width=0.40\textwidth]{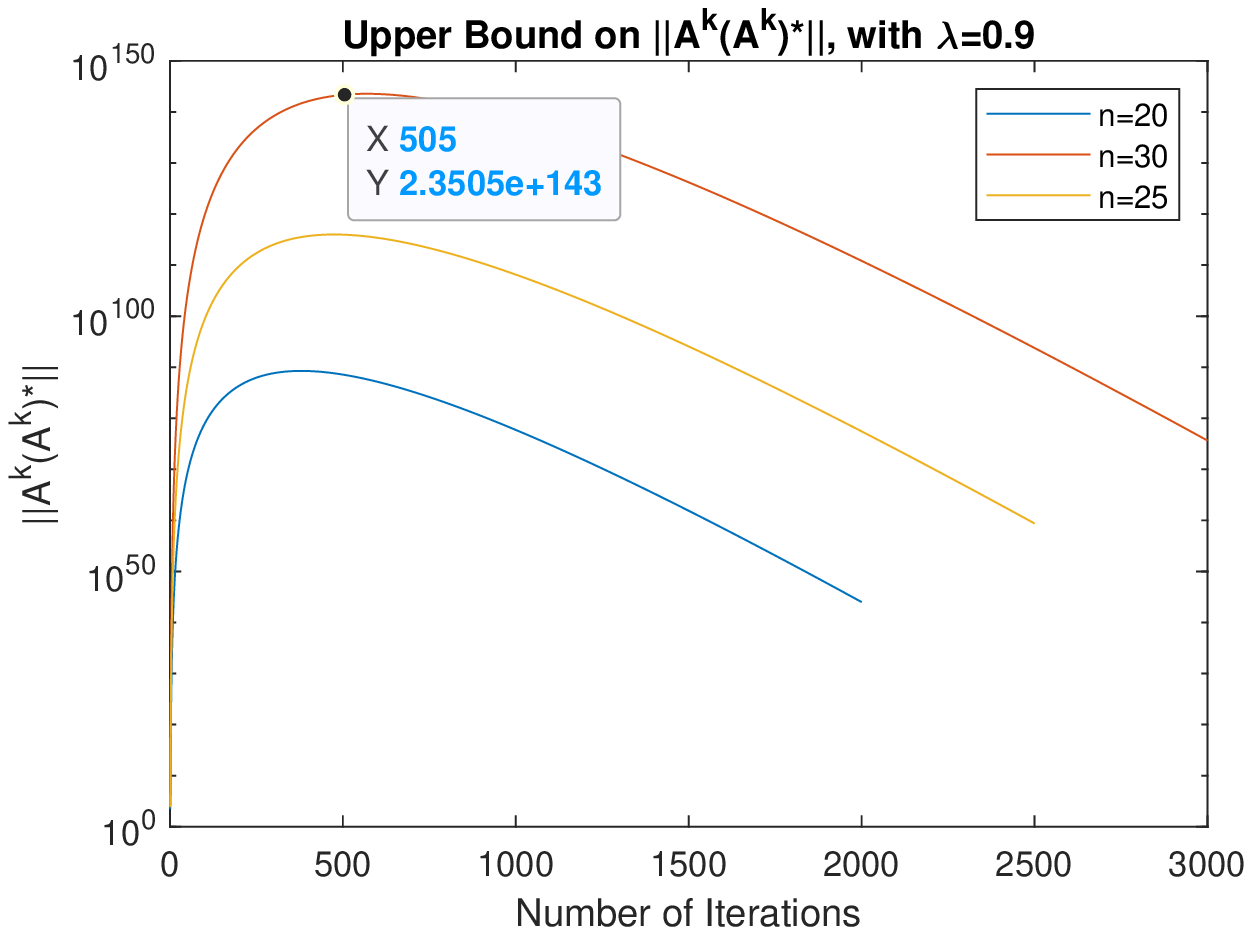}
    }
    \caption{(a) Setup: We fix the dimension of the state space to 50 and chose only one eigenvalue of 0.9, and observe the affect of decrease in spatial correlations by, first run: all the dynamics are represented by one Jordan block(invariant subspace), second run: divide them in two blocks and for the third run we divide dynamics into 3 different blocks, resulting in noticeable improvement in OLS error.(b) Setup: we fix two eigenvalues 0.9 and -0.75, with associated invariant subspaces of size n-3 and 3, respectively. Dimension of the underlying state space is varied from 5 to 70 and ratio of underlying state space and length of simulated trajectory is fixed at 1:3. Increasing the dimension of the state space, smallest singular value of stable non-diagonalizable dynamics scales better than that of a Gaussian ensemble(N x n rectangular matrix of i.i.d normals) (c)  but the largest singular value grows exponentially in dimension of the state space which explains error in OLS by using negative second moment identity (d) Stable n dimensional dynamical system with only one distinct eigenvalue and n-1 generalized eigenvectors, blow up of covariance.}   
    \label{fig:spatialcorr}
\end{figure}

\paragraph{Vulnerability of Stable, Non-Hermitian matrices to strong spatial correlations, leading to extremely strong temporal correlations:}
Position and magnitude alone of eigenvalues associated to a linear operator $A$, are not sufficient  to describe its' complete behavior (for the ease of exposition, throughout this paper we will assume that $A$ does not have any non-trivial null space). Knowledge about a linear operator $A$ is equivalent to knowing its' invariant subspaces (see e.g., \cite{kato2013perturbation}). Most of the following exposition is from \cite{naeem2023learning}; algebraic multiplicity of eigenvalues follow from determinant of the matirx. 
\begin{equation}
    \label{eq:det} det(zI-A)= \prod_{i=1}^{K} (z-\lambda_{i})^{m_i},
\end{equation}
where $\lambda_{i}$ are distinct with multiplicity $m_{i}$. Analysis becomes more involved when $dim[N(A-\lambda_{i}I)]<m_{i}$, which leads to invariants subspace (spanned by more that one linearly independent vector). Consequently, states space can be written as direct sum decomposition of $A-$ invariant subspaces 
\begin{equation}
\label{eq:directsumA}
    \mathbb{R}^{n}= M_{\phi(1)} \oplus M_{\phi(2)} \oplus \ldots \oplus M_{\phi(L)},
\end{equation}
 with respective orthogonal projetcions $[E_{\phi(i)}]_{i=1} ^{L}$ such that identity matrix can be written as:
\begin{equation}
    \label{eq:addId} I_{n}=E_{\phi(1)} \oplus E_{\phi(2)} \oplus \ldots \oplus E_{\phi(L)}
\end{equation} 
and $A$ as a composition of:
\begin{equation*}
    A=\bigoplus_{i=1}^{L} A_{\phi(i)}
\end{equation*}
where $\phi$ is a surjective map from $\{1,\ldots,L\}$ to $\sigma(A)$. $\phi$ is bijective iff eigenvectors span $\mathbb{R}^{n}$. \emph{In the case of gap between between algebraic and geometric multiplicity related to some element of $\sigma(A)$},
consider the invariant subspace $M_{\lambda}$, for some $\lambda \in \sigma(A)$, with algebraic multiplicity of $\lambda$ is $|B_{\lambda}|$ but only one linearly independent eigenvector $v_1$ such that $Av_1=\lambda v_1$. So we generate generalized eigenvector $v_2, v_3, \ldots, v_{|B_{\lambda}|}$ recursively as $(A-\lambda I)v_2=v_1$ and $(A-\lambda I)v_3=v_2$ and so on. We have the following $k$ -th step iteration:
\begin{align}
     & \nonumber A^{k}v_1 =\lambda^{k}v_1 \\ & \nonumber A^{k}v_2=\lambda^{k}v_{2} + \binom{k}{1}\lambda^{k-1}v_1 \\ & \nonumber A^{k}v_3=\lambda^{k}v_{3} + \binom{k}{1}\lambda^{k-1}v_2 +  \binom{k}{2} \lambda^{k-2}v_1 \\ & \nonumber \ldots= \ldots \\ & \label{eq:recursiveJordan}
    A^{k}v_{|B_{\lambda}|}= \lambda^{k}v_{|B_{\lambda}|}+ \binom{k}{1}\lambda^{k-1}v_{|B_{\lambda}|-1}+\ldots   + \binom{k}{|B_{\lambda}|-1}v_1,
\end{align}
which is conducive of strong spatial correlations for associated dynamics, also verified in the figure (a) of \ref{fig:spatialcorr}
\paragraph{Statistics of largest singular value:}Preceding argument hints at strong spatial correalation, combined with results from the section on diagonalizable dynamics one can guess singular values are distributed non-uniformly also shown in (b) and (c) of figure \ref{fig:spatialcorr}: least singular value scales better than $\sqrt{N}-\sqrt{n}$ but largest singular value is $O(e^{nN})$ that we aim to now formalize. 
\begin{theorem}  We can rigorously upper and lower bound norm of the $k-th$ iteration associated to action of matrix $A$ on invariant subspace $M_{\lambda}$, precisely given as:  
    \begin{align}
       \frac{1}{\sqrt{|B_{\lambda}|}}\min_{1 \leq i \leq |B_{\lambda}|} d(A_{M_{\lambda}}^{k} (i), A_{M_{\lambda}}^{k} (-i) ) \leq  \|A^{k}_{M_{\lambda}}\|_{2} \leq |\lambda|^{k} k^{|B_{\lambda}|}\sum_{m=0}^{|B_{\lambda}|-1} \frac{1}{|\lambda|^{m}} \overbrace{\leq}^{|\lambda| \in (0,1)} k^{|B_{\lambda}|}|B_{\lambda}|  |\lambda|^{k+1- |B_{\lambda}| } 
    \end{align}
where $A^{k}_{M_{\lambda}}:=A^{k}E_{\lambda}$  is the $k-$ th power of the Jordan block corresponding to invariant subspace $M_{\lambda}$, with algebraic representation $(A^{k}_{M_{\lambda}})_{ij}= \binom{k}{j-i} \lambda^{k-j+i}$. $A_{M_{\lambda}}^{k} (i)$ corresponds to the $i-th$ row of the $A^{k}_{M_{\lambda}}$ and $A_{M_{\lambda}}^{k} (-i)$ represents the span of all the rows except the $i-$ th one.    
\end{theorem}
\begin{proof}
First inequality follows from bounds on least singular value of $A^{k}_{M_{\lambda}}$, $\sigma_{|B_{\lambda}|} \big(A^{k}_{M_{\lambda}}\big)$ given in Lemma 1.3 of Chapter 1 in \cite{vu2014modern}:
\begin{equation}
    \frac{1}{\sqrt{|B_{\lambda}|}}\min_{1 \leq i \leq |B_{\lambda}|} d(A_{M_{\lambda}}^{k} (i), A_{M_{\lambda}}^{k} (-i) ) \leq \sigma_{|B_{\lambda}|} \big(A^{k}_{M_{\lambda}}\big) \leq \min_{1 \leq i \leq |B_{\lambda}|} d(A_{M_{\lambda}}^{k} (i), A_{M_{\lambda}}^{k} (-i) )
\end{equation}
General upper bound follows by a simple variation of bounds provided in \cite{tropp2001elementary} and last inequlaity under stability assumption is shown in \cite{naeem2023learning}.
\end{proof}
Now w.l.o.g assume that $|B_{\lambda}|=n$ for stable $\lambda$ , i.e., that is only one  eigenvalue with one eigenvector and $n-1$ generalized eigenvectors. Consequntly, $\|A^{k}_{M_{\lambda}}(A^{k}_{M_{\lambda}})^{T}\|_{2} \leq k^{2n} |\lambda|^{2k} n^2 \big(\frac{1}{|\lambda|}\big)^{2(n-1)}$ and for $k=O(n)$ grows as an exponential and then polynomially attaining a maxima $k=\Theta(\frac{n \ln n}{\ln(\frac{1}{|\lambda|})})$ and then starts decaying (see (d) of figure \ref{fig:spatialcorr}). Recall, conditioned on $x_0=0$, we can express $N-th$ realization of the signal as: $x_{N}= \sum_{k=1}^{N} A_{M_{\lambda}}^{N-k} \eta_{k-1}$, along with preceding analysis explains blow up of largest singular value (see (c) of figure \ref{fig:spatialcorr}), which implies that distance between some or many of the row vectors and conjugate  hyperplanes (defined as in \ref{thm:neg_2nd_moment_ineq}) is becoming small.

\section{Conclusion and future work:}
\label{sec:conc}
In this paper we proposed a novel analysis for system identification, via running a single trajectory. Our analysis is based on concentration of measure phenomenon and tools developed in random matrix theory, making it naturally apt for high dimensional dynamical systems. Some of the important results and suggestions: qualitatively different behavior of diagonalizable and non-hermitian dynamical systems. Failure of OLS on non-hermitian systems, stems from the strong spatial correlations. Explosive diagonalizable dynamics, have large variance and distance between associated hyperplane and random vector can be small with high probability, thus deteriorating performance of OLS. Our analysis offers new and deep insights into learning for dynamical systems by studying the problem in Frobenius norm. However, more effort is required to bring our analysis on non-hermitian systems at the same theoretical footing as in diagonalizable case, which is currently a work in progress. 

\bibliography{neurips_2023}
\bibliographystyle{plainnat}

\end{document}